\documentclass[12pt]{article}
\usepackage{graphicx} % Required for inserting images
\usepackage{amsmath,amssymb,enumerate,xcolor}
\usepackage[hidelinks]{hyperref}

\title{Self-adjoint quantization of Stäckel integrable systems.}
\author{ Jonathan M.\ Kress\footnote{School of Mathematics and Statistics, UNSW Sydney, NSW 2052, Australia}  \quad \& \quad  Vladimir S.\ Matveev\footnote{
Institut f\"ur Mathematik, Friedrich Schiller Universit\"at Jena,
07737 Jena,  Germany  \ \ \quad {\tt  vladimir.matveev@uni-jena.de}}}

\newtheorem{theorem}{Theorem}
\newtheorem{remark}[theorem]{Remark}

\newcommand{\weg}[1]{}

\begin{document}

\maketitle
\begin{abstract} 
We show that  quadratic Hamiltonians in involution  coming from a Stäckel system are  quantizable, in the sense that one can construct commutative self-adjoint 
operators whose symbols are the quadratic Hamiltonians. Moreover, they allow multiplicative separation of variables.  This proves a conjecture explicitly formulated in \cite{BMD16}.

\end{abstract}

\section{Introduction}
We work in a simply-connected neighbourhood   $W\subset  \mathbb{R}^n$ with coordinates $x^1,...,x^n$ and denote the corresponding momenta by $p_1, \ldots , p_n$.

It is known that given a nondegenerate $n\times n$ matrix $S$ such that for every $i$ the components $S_{ij}$ in the $i$ row are functions of $x^i$ only, one can construct a family of functions $H_1,\ldots, H_n: T^*W\to \mathbb{R}$, given by (\ref{eq:St_intro}), which commute with respect to the standard Poisson bracket 
$$\{H, F\}:=\sum^n_{i=1}\left( \frac{\partial H}{\partial p_i} \frac{\partial F}{\partial x^i} -\frac{\partial H}{\partial x^i} \frac{\partial F}{\partial p_i} \right).$$   This   construction,  known as the   {\it St\"ackel construction},  appeared already   in \cite[\S\S 13-14]{Liouville}, see also discussion in \cite[pp. 703--705]{Luetzen} and was later independently found   and intensively studied by P.
Stäckel in e.g.   \cite{disser}, and  further  by e.g. L.P.  Eisenhart \cite{eisenhart}.

Let us recall the construction. Consider the functions $H_\alpha$, $\alpha =1,\ldots, n$,  given by the following system of linear equations
 \begin{equation} 
 \label{eq:St_intro}
S\mathbb{ H}  = \mathbb{P},  
\end{equation} 
where  $\mathbb{H}= \left( H_1, H_2,\ldots ,H_n\right)^\top$ and $\mathbb{P}= \left(p_1^2, p_2^2,\ldots ,p_n^2\right)^\top$.  It is known that the functions $H_\alpha$ are in involution, that is, $\{H_\alpha, H_\beta\}=0 $, and that they are functionally independent almost everywhere, so they generate  an integrable system.  We will call them  {\it Hamiltonians}.  They  are clearly quadratic in the momenta $p$.

By    {\it quantization} of the  integrable system generated by $H_1, H_2,\ldots ,H_n $ we understand a system of commutative   second order differential operators $\hat H_1, \hat H_2, \ldots , \hat H_n$ such that for every $i=1,\ldots, n$ the symbol of $\hat H_i$ is $H_i$. We say that the quantization  is {\it self-adjoint}, if all the
operators $\hat H_i$ are self-adjoint  with respect to the following inner product  on the space of 
compactly  supported smooth functions on $W$: 
\begin{equation} \label{eq:2}
\langle f, h \rangle_{\phi}= \int_W f h \phi(x) dx^1\wedge \ldots  \wedge dx^n.
\end{equation}
Here $\phi(x)dx^1\wedge \ldots  \wedge dx^n$ should be viewed as   the  volume form which defines the inner product \eqref{eq:2}. We assume that 
$\phi$  never vanishes. 

Using integration by parts (or Stokes' theorem) one immediately sees that 
a self-adjoint  with respect to \eqref{eq:2}  second order operator  with symbol $K^{ij}p_ip_j$ is necessarily of the form  
\begin{equation} \label{eq:3}
\sum_{i,j}\frac{1}{\phi} \frac{\partial }{\partial x^i} K^{ij} \phi \frac{\partial }{\partial x^j} + U, 
\end{equation}
where $U:W\to \mathbb{R}.$

\begin{remark} {\phantom{ 1}}
\begin{enumerate}\item The  form   \eqref{eq:3} is covariant providing $\phi$ changes as the  coefficient of an  $n$-form  $\phi(x)dx^1\wedge \ldots  \wedge dx^n$ and $V$ as a function. 
   \item An equivalent way to write the first term in the  formula \eqref{eq:3} is 
\begin{equation} \label{eq:3b}
\sum_{i,j}\tilde \nabla_i K^{ij}  \tilde \nabla_j, 
\end{equation}
where $\tilde \nabla$ is any symmetric affine connection such that the $n$-form  $\phi(x)dx^1\wedge \ldots  \wedge dx^n$ is parallel. 
\item  If $K^{ij}=g^{ij}$ for a Riemannian or pseudo-Riemannian metric $g_{ij}$  and   $\tilde \nabla$ is the   Levi-Civita connection of $g$, 
the formula \eqref{eq:3b} gives us the Laplace-Beltrami operator of $g$. 
\end{enumerate}
\end{remark}

  We denote   the symmetric $(2,0)$-tensor corresponding to the functions $H_\alpha$   by $H_{(\alpha)}^{ij}$. In the coordinates $x^i$, $H_{(\alpha)}^{ij}$
  is given by the diagonal matrix and its $i$th diagonal component  $H_{(\alpha)}^{ii}$
equals  $\tfrac{\Delta_{\alpha i}}{\det(S)}$, where we denote by $\Delta_{ij} $ the components of the comatrix  of $S$.

Our first  result is the  following theorem: 

\begin{theorem} \label{thm:main}
For the  Stäckel integrable system  constructed by an arbitrary  Stäckel matrix $S$    and for  the function  $\phi:= \det(S)$, 
the second order differential operators \begin{equation} \label{eq:last}
\hat H_\alpha = \sum_{i,j}\frac{1}{\phi} \frac{\partial }{\partial x^i} H_{(\alpha)}^{ij} \phi \frac{\partial }{\partial x^j}, 
\end{equation} commute. 
\end{theorem}

\begin{remark}
The choice of $\phi$ is not unique as a trivial change of the separable coordinates associated with the St\"ackel system, i.\ e.\ $\hat x^i=f_i(x^i)$, preserves the St\"ackel property and gives $d\hat x^i=f'_i(x^i)dx^i$, for $i=1,\ldots,n$, and hence multiplies $\phi$ by a product of one-variable functions.    
\end{remark}

Clearly, the symbols of the operators  $\hat H_\alpha$ are the Hamiltonians  $H_\alpha$.

Special cases of Theorem  \ref{thm:main} were known before.
The most known   one  is when  
the matrix $S$ is the Vandermonde matrix
\begin{equation}\label{eq:dim2}
S= \begin{pmatrix}(f_1(x^1))^{n-1}  &(f_1(x^1))^{n-2} & \cdots & 1 \\
(f_2(x^2))^{n-1}  &(f_2(x^2))^{n-2} & \cdots & 1   \\ 
\vdots &     \vdots        &        & \vdots \\ 
(f_n(x^n))^{n-1}  &(f_n(x^n))^{n-2} & \cdots & 1. 
\end{pmatrix}
\end{equation}
This case was done e.g. in \cite{M2020, MT01}. In this case, the function $\phi=\det(S)$ coincides, up to a sign, with $\sqrt{|\det g|}$, where $g$ is the metric such that 
$\sum_{i,j} g^{ij}p_ip_j= H_1$. The  condition that the   matrix $S$ is nondegenerate implies here and also in the next case 
that at no point $f_i(x^i)= f_j(x^j)$.   Note that many famous integrable systems such that geodesic flow of the ellipsoid or of the Poisson sphere can be obtained with the help of Stäckel matrix $S$ given by \eqref{eq:dim2}, see e.g. \cite{MT01}. 
We remark here that  in dimension $n=2$, any Stäckel system can be obtained using
$S$ given by \eqref{eq:dim2}, so the two-dimensional  case is also known.

Another known special case corresponds to the Stäckel matrix 
$$
S= \begin{pmatrix}(f_1(x^1))^{\gamma_{n-1}}  &(f_1(x^1))^{\gamma_{n-2}} & \cdots & (f_1(x^1))^{\gamma_{0}} \\
(f_2(x^2))^{\gamma_{n-1}}  &(f_2(x^2))^{\gamma_{n-2}} & \cdots & (f_2(x^2))^{\gamma_{0}}    \\ 
\vdots &     \vdots        &        & \vdots \\ 
(f_n(x^n))^{\gamma_{n-1}}  &(f_n(x^n))^{\gamma_{n-2}} & \cdots & (f_n(x^n))^{\gamma_{0}}  
\end{pmatrix},
$$
where $\gamma_i$ are mutually different integer numbers. The existence of  a self-adjoint  quantization  (and of  multiplicative separability, see Theorem \ref{thm:3} below) for the corresponding integrable system  
was established  in \cite{BMD16}, see also \cite[\S 8.2.2]{Blaszakbook}.

One more special case is related to orthogonal separation of varables  in spaces of constant curvature, that is, when the Stäckel matrix  is such that 
one of the Hamiltonians  $H_\alpha$, say $H_1$, is the  kinetic energy corresponding to a metric $g$ of constant sectional curvature.   Recall  for  a  metrics constant curvature, the Robertson condition is trivially fulfilled, so 
the quantization  corresponding to $\phi= \sqrt{|\det(g)|}$ produces commutative operators, see e.g. \cite{BCR}.     The corresponding Stäckel matrices  were described in \cite[Theorem 1.4]{BKM22}. 

Clearly, in all the cases above, elementary transformations with the rows of the Stäckel matrix correspond to taking linear combinations of the corresponding integrals and do  not affect quantization. 

We did not find other special cases in the literature.

Let us now discuss the possibility of adding potentials to the ``pure kinetic'' Hamiltonians  $H_\alpha$, such that the integrability, and the corresponding quantization, is preserved.  Recall that, in the coordinates $x^i$ from above, it is well-known what functions  $U_i:U\to \mathbb{R}$ have the property that 
the functions $H_i+ U_i$, with $H_i$ as  above,  are in involution.  The freedom in constructing such  functions $U_i$ is the choice of $n$ functions $V_1(x^1), V_2(x^2),\ldots, V_n(x^n)$ of the indicated variables, and the functions $U_i$ are given by the formula 
\begin{equation} 
 \label{eq:St_intro1}
S\mathbb{ U}  = \mathbb{V},  
\end{equation} 
where  $\mathbb{U}= \left( U_1,U_2,\ldots ,U_n\right)^\top$ and $\mathbb{V}= \left(V_1(x^1), V_2(x^2),\ldots ,V_n(x^n)\right)^\top$.

\begin{theorem} \label{thm:2}
For any Stäckel integrable system constructed from the Stäckel matrix $S$, the differential operators 
$\check H_i= \hat H_i + U_i$,  for $i=1,\ldots n$, commute.   Moreover, if the operators $\check H_i= \hat H_i + U_i$ for some functions $U_i$ commute, then the functions $U_i$ are constructed by certain functions $V_1(x^1), V_2(x^2),\ldots, V_n(x^n)$ by \eqref{eq:St_intro1}. 
\end{theorem}
   As it will be clear from the proofs  of Theorems \ref{thm:main} and \ref{thm:2},  in the coordinates $x^1,\ldots, x^n$, the operators  $\check H_i$ are given by  
    \begin{equation} \check H_{\alpha} = \sum_{i} H^{ii}_{(\alpha)} \left[ \left(\frac{\partial}{\partial x^i}\right)^2 + V_i(x^i) \right].
\label{HwithV}
    \end{equation}

Let us now discuss (multiplicative, or quantum in the terminology of \cite{BMD16}) separability of the corresponding joint eigenfunction  problem. We take arbitrary $E_1,\ldots, E_n\in \mathbb{R}$ (or $\mathbb{C}$) and  consider the  following PDE-system on the unknown function $\Psi= \Psi(x^1,\ldots, x^n)$: 
\begin{equation}\label{eq:sys1} \begin{array}{rcl}
    \check H_1\Psi &=&E_1\Psi ,  \\ & \vdots & \\   \check H_n\Psi &=&E_n\Psi .
\end{array}\end{equation}

\begin{theorem} \label{thm:3}
  Consider the functions $\psi_1(x^1), \psi_2(x^2), \ldots , \psi_n(x^n)$ of the indicated variables   such that
  for $\alpha=1,\ldots, n$  the function $\psi_\alpha(x^\alpha)$
  satisfies  the following ODE:
\[
\left(S_{\alpha 1}E_1 + S_{\alpha 2}E_2 + \cdots + S_{\alpha n}E_n \right)\psi_\alpha(x^\alpha)
 =  \left[\left(\frac{d}{d x^\alpha}\right)^2 + V_\alpha(x^\alpha)\right]\psi_\alpha(x^\alpha).
\]
Then, the function $$\Psi(\mathbf x) = \prod_{\alpha=1}^n\psi_\alpha(x^\alpha)$$ 
satisfies \eqref{eq:sys1}. 
\end{theorem}
Theorem \ref{thm:3} proves  the conjecture explicitly formulated in \cite[ page 461]{BMD16}.

 \subsection*{Acknowledgement. } V.\ M.\ thanks the DFG (projects 455806247 and 529233771), and the ARC  (Discovery Programme DP210100951) for their support. J.\ K. thanks V.\ M.\  for hospitality during a visit to Friedrich Schiller Universit\"at Jena, where most of this work was carried out, and Friedrich Schiller Universit\"at Jena for partial financial support.

\section{ Proof of Theorem \ref{thm:main}.}
 Clearly, $\Delta_{\alpha i} $  does not depend on  the variable $x^i$.
The corresponding operator   
$\hat H_\alpha$ is then given by 
\begin{equation}\label{eq:4}
\hat H_\alpha=  \sum_i \frac{1}{\phi} \frac{\partial }{\partial x^i} \frac{\Delta_{\alpha i}}{\phi}  \phi \frac{\partial }{\partial x^i}        =   \sum_i \frac{\Delta_{\alpha i}}{\phi}  \left(\frac{\partial  }{\partial  x^i} \right)^2 =\sum_i  H_{(\alpha)}^{ii}  \left(\frac{\partial   }{\partial  x^i } \right) ^2.        
\end{equation}
Let us now commute two such operators, $\hat H_\alpha$ and $\hat H_\beta$. The result clearly has no   terms of the form 
$$ \textrm{function} \cdot  \frac{\partial   }{\partial  x^i }.   $$ 
As $\{H_\alpha, H_\beta\}=0,  $ the result has no  terms of the form  (we call them {\it third  order terms})
$$ \textrm{function} \cdot  \frac{\partial^3    }{\partial  x^i \partial  x^j  \partial  x^k }.   $$ 
From the other side,  in view of \eqref{eq:4},  the sum of the third 
order terms of  the commutator $[\hat H_\alpha, \hat H_\beta]$  
is 
\begin{equation} \label{eq:5}  
\frac{1}{\phi}    \sum_{i,j} \left( \frac{  \partial}{\partial  x^i}  \frac{\Delta_{\alpha i} \Delta_{\beta j} -\Delta_{\alpha j} \Delta_{\beta i} }{\phi}\right) 
\frac{\partial^3 }{\partial x^i \partial x^j \partial x^j } .      
\end{equation}
This implies that  for  every  $i,j$  \begin{equation} \label{eq:6} 
\frac{  \partial}{\partial  x^i}  \frac{\Delta_{\alpha i} \Delta_{\beta j} -\Delta_{\alpha j} \Delta_{\beta i} }{\phi}  =0 .\end{equation}

 Thus, the commutator is given by 
 $$
 [\hat H_\alpha, \hat H_\beta]=\frac{1}{\phi} \sum_j \left( \sum_i \left(\frac{  \partial}{\partial  x^i}\right)^2 \frac{\Delta_{\alpha i} \Delta_{\beta j} -\Delta_{\alpha j} \Delta_{\beta i} }{\phi}\right)   \left(\frac{\partial }{  \partial  x^j}\right)^2.  
 $$
 Its coefficients at $\left(\tfrac{\partial }{  \partial  x^j}\right)^2$ is therefore  the sum of $\tfrac{\partial }{\partial x^i}-$derivatives of the left hand sides of 
 \eqref{eq:6} and vanishes identically, implying  $[\hat H_\alpha, \hat H_\beta]=0$ as we claimed. Theorem \ref{thm:main} is proved.

\section{ Proof of Theorem \ref{thm:2}.}

We wish to show that $\check H_{(\alpha)}=\hat H_{(\alpha)} + U_\alpha$ and $\check H_{(\beta)}=\hat H_{(\beta)}+U_\beta$, given by \eqref{HwithV},
commute. Since $[\hat H_{(\alpha)},\hat H_{(\beta)}]=0$ by Theorem \ref{thm:main} and scalar functions commute,
\[
[\check H_{(\alpha)},\check H_{(\beta)}]
= [\hat H_{(\alpha)},U_\beta]  - [\hat H_{(\beta)},U_\alpha].
\]
The highest (second) derivative terms will clearly cancel.  To check that all other terms vanish, we write them in the coordinates $x^i$ associated with the St\"ackel system as,
\begin{eqnarray*}
    [\check H_{(\alpha)},\check H_{(\beta)}]
    &=& \left[ \sum_{j=1}^n  \frac{\Delta_{\alpha j}}
    {\phi } \left(\frac{\partial}{\partial x^j}\right)^2 , \sum_{i=1}^n \frac{\Delta_{\beta i}}
    {\phi } V_i(x^i) \right] 
%    \\
%    & & {} \ 
    % - \left[ \sum_{j=1}^n  \frac{\Delta_{\beta j}}
    % {\phi } \left(\frac{\partial}{\partial x^j}\right)^2, \sum_{i=1}^n \frac{\Delta_{\alpha i}}
    % {\phi } V_i(x^i) \right] \\
    - (\alpha \longleftrightarrow \beta) \\
    &=&  \sum_{j=1}^n \sum_{i=1}^n \left(
      \frac{2\Delta_{\alpha j}}{\phi } 
      \frac{\partial}{\partial x^j}  \left(\frac{\Delta_{\beta i}}{\phi}\right) V_i(x^i)\frac{\partial}{\partial x^j}
     + \frac{2\Delta_{\alpha j}\Delta_{\beta i}}
    {\phi^2 } \delta_{ij}V'_i(x^i)\frac{\partial}{\partial x^j}  \right.\\
    & & {}
    + \frac{\Delta_{\alpha j}}
    {\phi } \left(\frac{\partial}{\partial x^j}\right)^2\left( \frac{\Delta_{\beta i}}
    {\phi }\right) V_i(x^i)
    + 2\frac{\Delta_{\alpha j}}
    {\phi } \frac{\partial}{\partial x^j}\left( \frac{\Delta_{\beta i}}
    {\phi }\right) \delta_{ij}V'_i(x^i) \\
   & & {} \ \left.   + \frac{\Delta_{\alpha j}\Delta_{\beta i}}
    {\phi^2 } \delta_{ij} V''_i(x^i)
 \phantom{\frac{\partial}{\partial x^j}} - (\alpha \longleftrightarrow \beta )  \right)
\end{eqnarray*}
The terms with a factor of $\delta_{ij}$ are symmetric in $\alpha$ and $\beta$, at least after summing over $i$, and so cancel.  Noting that $\Delta_{\alpha i}$ does not depend on $x^i$, the remaining terms can be collected as,
\[
\sum_{j=1}^n \sum_{i=1}^n \left(
\frac{2V_i(x^i)}\phi \frac{\partial}{\partial x^j}\left(\frac{\Delta_{\alpha j}\Delta_{\beta i} - \Delta_{\beta j}\Delta_{\alpha i}}\phi\right)\frac{\partial}{\partial x^j}
+ \frac{V_i(x^i)}\phi\left(\frac{\partial}{\partial x^j}\right)^2\left(\frac{\Delta_{\alpha j}\Delta_{\beta i} - \Delta_{\beta j}\Delta_{\alpha i}}\phi\right)
\right),
\]
and these vanish by virtue of (\ref{eq:6}).  Hence $[\check H_{(\alpha)},\check H_{(\beta)}] = 0$.

Let us now prove Theorem in the other direction.  We will first consider the conditions for two Hamiltonians 
$$
H_\alpha+ U_\alpha = \sum_i H_{(\alpha)}^{ii} p_i^2 + U_\alpha \ \textrm{and} \ H_\beta+ U_\beta = \sum_i H_{(\beta)}{}_{ii} p_i^2 + U_\beta 
$$
to be in involution. We assume that  the quadratic parts $H_\alpha, H_\beta$ came from the Stäckel construction and $U_\alpha, U_\beta$ are unknown functions on $W$.  Since $\{H_\alpha,H_\beta\}=0 =\{U_\alpha, U_\beta\},$ we obtain 
\begin{eqnarray*}
\{H_{\alpha }   + U_\alpha,  H_{\beta }  + U_\beta \} &= &  \{H_\alpha, U_\beta\}-\{H_\beta, U_\alpha\} \\ 
  & =&   2  \sum_{s } \left(H_{(\alpha)}^{ss} p_s  \frac{\partial U_\beta}{\partial x^s}   - H_{(\beta)}^{ss} p_s  \frac{\partial U_\alpha}{\partial x^s}   \right).  
\end{eqnarray*}
We see that the Hamiltonians are in involution,  if     for any $s  =1,\ldots, n$,  the ``Benenti'' condition 
\begin{equation}  \label{eq:benenti}
H_{(\alpha)}^{ss}   \frac{\partial U_\beta}{\partial x^s}   - H_{(\beta)}^{ss}   \frac{\partial U_\alpha}{\partial x^s}=0\end{equation} is fulfilled. 
Next, consider the commutator $[\check H_\alpha, \check H_\beta]$, where again $U_\alpha, U_\beta$ are arbitrary functions. In view of   $[\hat H_\alpha,\hat H_\beta]=0 =[U_\alpha, U_\beta]$,  we have 
\begin{eqnarray*}
[\hat H_{\alpha }   + U_\alpha,  \hat H_{\beta }  + U_\beta ]&= &  [\hat H_\alpha, U_\beta]+[U_\alpha, \hat H_\beta] \\ 
  & =&      \sum_{s } \left(H_{(\alpha)}^{ss}   \frac{\partial U_\beta}{\partial x^s}  - H_{(\beta)}^{ss}  \frac{\partial U_\alpha}{\partial x^s}  \right) \tfrac{ \partial }{\partial x^s }  \\ 
  &+&    \sum_{s } \left(H_{(\alpha)}^{ss}   \frac{\partial^2 U_\beta}{\partial (x^s)^2} - H_{(\beta)}^{ss}  \frac{\partial^2 U_\alpha}{\partial (x^s)^2}\  \right). 
\end{eqnarray*}
We see that the condition $[\check H_\alpha, \check H_\beta]=0$ implies  \eqref{eq:benenti}. Thus,  if the  operators commute, the functions $H_{\alpha }   + U_\alpha$ and $  H_{\beta }  + U_\beta$ are in involution implying our  claim. 
  
\section{ Proof of Theorem \ref{thm:3}.}

We seek  solutions $\Psi$ to the system \eqref{eq:sys1}. 
The standard trick, see e.g. \cite[proof of Theorem 8]{BMD16},  based on the usual combination of the Hamiltonians leads to uncoupled equations.  That is, using (\ref{HwithV}),
\begin{eqnarray*}
    \lefteqn{\left(S_{\alpha 1}\check H_1 + S_{\alpha 2}\check H_2 + \cdots + S_{\alpha n}\check H_n \right)\Psi(x^1,x^2,\ldots x^n)} \\
    % &=& \sum_{j=1}^n S_{\alpha j}\sum_{i=1}^n H_{(j)}^{ii}\left[\left(\frac{\partial}{\partial x^i}\right)^2 + V_i(x^i)\right]\Psi(x^1,x^2,\ldots x^n) \\
    &=& \sum_{j=1}^n S_{\alpha j}\sum_{i=1}^n \frac{\Delta_{j i}}{\phi}\left[\left(\frac{\partial}{\partial x^i}\right)^2 + V_i(x^i)\right]\Psi(x^1,x^2,\ldots x^n) \\
    &=& \sum_{i=1}^n \sum_{j=1}^n S_{\alpha j}\frac{\Delta_{j i}}{\phi}\left[\left(\frac{\partial}{\partial x^i}\right)^2 + V_i(x^i)\right]\Psi(x^1,x^2,\ldots x^n) \\
    &=& \sum_{i=1}^n \delta_{\alpha i}\left[\left(\frac{\partial}{\partial x^i}\right)^2 + V_i(x^i)\right]\Psi(x^1,x^2,\ldots x^n) \\
    &=& \left[\left(\frac{\partial}{\partial x^\alpha}\right)^2 + V_\alpha(x^\alpha)\right]\Psi(x^1,x^2,\ldots x^n)
\end{eqnarray*}
Now, with the Ansatz, $\Psi(x^1,x^2,\ldots,x^n)=\psi_1(x^1)\psi_2(x^2)\ldots\psi_n(x^n)$, we have the $n$ separated ordinary differential  equations,
\[
\left(S_{\alpha 1}E_1 + S_{\alpha 2}E_2 + \cdots + S_{\alpha n}E_n \right)\psi_\alpha(x^\alpha)
 =  \left[\left(\frac{d}{d x^\alpha}\right)^2 + V_\alpha(x^\alpha)\right]\psi_\alpha(x^\alpha),
\]
for $\alpha=1\ldots n$. This system of (uncoupled)  ODEs on the functions $\phi_\alpha(x^\alpha)$ is clearly equivalent to the system  \eqref{eq:sys1}. 

% \begin{thebibliography}{99}

% \end{thebibliography}
\bibliographystyle{plain}
% \bibliography{quantisation_staeckel_systems}

\begin{thebibliography}{10}

\bibitem{BCR}
S.~Benenti, C.~Chanu, and G.~Rastelli.
\newblock Remarks on the connection between the additive separation of the {H}amilton-{J}acobi equation and the multiplicative separation of the {S}chr\"{o}dinger equation. {I}. {T}he completeness and {R}obertson conditions.
\newblock {\em J. Math. Phys.}, 43(11):5183--5222, 2002.

\bibitem{Blaszakbook}
Maciej B{\l}aszak.
\newblock {\em Quantum versus classical mechanics and integrability problems---towards a unification of approaches and tools}.
\newblock Springer, Cham, 2019.

\bibitem{BMD16}
Maciej B{\l}aszak, Krzysztof Marciniak, and Ziemowit Doma\'nski.
\newblock Separable quantizations of {S}t\"ackel systems.
\newblock {\em Ann. Physics}, 371:460--477, 2016.

\bibitem{BKM22}
Alexey~V. Bolsinov, Andrey~Yu. Konyaev, and Vladimir~S. Matveev.
\newblock Orthogonal separation of variables for spaces of constant curvature.
\newblock {\em Forum Math.}, 37(1):13--41, 2025.

\bibitem{eisenhart}
L.~P. Eisenhart.
\newblock Separable systems of {S}tackel.
\newblock {\em Ann. of Math. (2)}, 35(2):284--305, 1934.

\bibitem{Liouville}
J.~Liouville.
\newblock M\'emoire sur l'int\'egration des \'equations diff\'erentielles du mouvement d'un nombre quelconque de points mat\'eriels.
\newblock {\em Journal de Math\'ematiques Pures et Appliqu\'ees}, 1e s{\'e}rie, 14:257--299, 1849.

\bibitem{Luetzen}
J.~L\"{u}tzen.
\newblock {\em Joseph {L}iouville 1809--1882: master of pure and applied mathematics}, volume~15 of {\em Studies in the History of Mathematics and Physical Sciences}.
\newblock Springer-Verlag, New York, 1990.

\bibitem{M2020}
V.~S. Matveev.
\newblock Quantum integrability for the {B}eltrami-{L}aplace operators of projectively equivalent metrics of arbitrary signatures.
\newblock {\em Chebyshevski\u i\ Sb.}, 21(2):275--289, 2020.

\bibitem{MT01}
Vladimir~S. Matveev and Peter~J. Topalov.
\newblock Quantum integrability of {B}eltrami-{L}aplace operator as geodesic equivalence.
\newblock {\em Math. Z.}, 238(4):833--866, 2001.

\bibitem{disser}
P.~St\"ackel.
\newblock Die integration der hamilton-jacobischen differentialgleichung mittelst separation der variablen.
\newblock {\em Habilitationsschrift, Universit\"at Halle}, 1891.

\end{thebibliography}

\end{document}